\documentclass[11pt]{amsart}
\usepackage[utf8]{inputenc}
\usepackage[T1]{fontenc}
\usepackage{amsmath,amssymb,amsthm,mathrsfs}
\usepackage{booktabs}
\usepackage[margin=2.8cm]{geometry}
\usepackage[colorlinks=true,linkcolor=blue,citecolor=blue,urlcolor=blue]{hyperref}

\newtheorem{theorem}{Theorem}[section]
\newtheorem{proposition}[theorem]{Proposition}
\newtheorem{lemma}[theorem]{Lemma}
\newtheorem{corollary}[theorem]{Corollary}

\theoremstyle{definition}
\newtheorem{definition}[theorem]{Definition}
\newtheorem{example}[theorem]{Example}
\newtheorem{remark}[theorem]{Remark}
\newtheorem{question}[theorem]{Question}

\newcommand{\Pc}{\mathcal P}
\newcommand{\Hc}{\mathcal H}
\newcommand{\Lc}{\mathcal L}
\newcommand{\Uc}{\mathcal U}
\newcommand{\Ind}{\operatorname{Ind}_p}
\newcommand{\Ip}{I_p}
\newcommand{\Ipd}{I_{pd}}
\newcommand{\dd}{\,\mathrm{d}}
\newcommand{\Var}{\operatorname{Var}}
\newcommand{\h}{\mathfrak h}

\title[Index of positivity of linear functionals]{On the index of positivity of linear functionals:\\ a complement to the Hermite case and new examples}

\author{Wathek Chammam}
\address{Department of Mathematics, College of Science, Majmaah University, Al-Majmaah 11952, Saudi Arabia}
\email{w.chammam@mu.edu.sa}


\begin{document}

\begin{abstract}
We prove that the Dirac functional $\delta$ is the only normalized linear functional on the space of real polynomials whose index of positivity is infinite. This complements the recent study of the Hermite case by R.~Sfaxi and F.~Marcell\'an [Bull. Malays. Math. Sci. Soc. 48 (2025) 138]: the index of the Hermite functional $\Hc$ is finite, and so are those of $\Lc(-\tfrac12)$, $e^{\Hc}$ and $e^{\Lc(-1/2)}$, for which the value $+\infty$ was stated there. The key observation is that real Cauchy powers act on log-moments by scaling. This yields explicit upper bounds for the index in terms of the first four moments. In particular $\Hc^n$ is not even positive for any integer $n\ge 3$, and $2\le \Ind(\Hc)\le 5$. We identify the point of the original argument that needs to be modified, a free integration constant in the construction of the weights of $\Hc^n$, and we give the adjusted statements. Finally, we compute the positivity sets exactly for two new families: the real powers of $\delta_c$, which are Jacobi functionals on $[0,c]$, and the real powers of $\delta_c\delta_{-c}$, which are generalized Gegenbauer functionals containing the Chebyshev functional $\mathcal T$. For the latter, $\Ind(\mathcal T)=2$. These examples show that our bounds are sharp and that every value in $[0,+\infty]$ is the index of some functional.
\end{abstract}

\maketitle

\section{Introduction}

Let $\Pc=\mathbb R[x]$ and let $\Pc'$ be its topological dual, which coincides with the algebraic dual. For $u\in\Pc'$ we write $(u)_n=\langle u,x^n\rangle$, and $u$ is \emph{normalized} if $(u)_0=1$. We use the notation of \cite{SM}: $uv$ is the Cauchy product, $e^u$ the Cauchy exponential of \cite{MS23}, $\delta_c$ the Dirac functional at $c$ and $\delta=\delta_0$. For the general theory of linear functionals and orthogonal polynomials we refer to \cite{Ch,GMM,Ma}. For a normalized $w$ there is a unique $v\in\Pc'$ with $(v)_0=0$ and $w=e^v$, and the real power is $w^\lambda:=e^{\lambda v}$. Following \cite{S24,SM}, set
\[
\Ip(w)=\{\lambda\ge 0:\ w^\lambda \text{ is positive}\},\qquad
\Ipd(w)=\{\lambda> 0:\ w^\lambda \text{ is positive-definite}\},\qquad
\Ind(w)=\sup \Ip(w).
\]
Since $w^0=\delta$, we always have $0\in \Ip(w)$, and $\Ind(\delta)=+\infty$. The index of positivity was introduced in \cite{S24}, where the question was raised whether there exist other functionals with infinite index. The paper \cite{SM} develops the tools needed to study this question, in particular the continuity of the Cauchy exponential, its series expansion and the differential relations satisfied by the powers of the Hermite functional $\Hc$, and it proposes an affirmative answer. Its Theorem~4.7 states that $\Hc^n$ is positive-definite for every integer $n\ge1$, whence $\Ind(\Hc)=+\infty$. From this it deduces $\Ind(\Lc(-\tfrac12))=+\infty$ (Section~4.3), and, through \cite[Theorem~3.5]{SM}, $\Ind(e^{\Hc})=\Ind(e^{\Lc(-1/2)})=+\infty$ (Corollaries~4.8 and 4.9).

The present note builds on the framework and on the tools of \cite{S24,SM}. Its starting point is the observation that real Cauchy powers act on log-moments by scaling, which leads to a different conclusion: the index is finite for every normalized functional other than $\delta$. The statements quoted above therefore need to be modified. Our main results are the following.

\begin{enumerate}
\item \textbf{Main theorem} (Section~\ref{sec:main}). For every normalized $w\neq\delta$ we have $\Ind(w)<+\infty$, with an explicit bound in terms of $(w)_1,\dots,(w)_4$ (Theorem~\ref{thm:main}). Hence the answer to the question of \cite{S24} is negative.
\item \textbf{The Hermite and Laguerre cases} (Section~\ref{sec:counter}). For the nonnegative polynomial $(x^2-3)^2$ one has $\langle \Hc^6,(x^2-3)^2\rangle=-\tfrac34$. Exact Hankel determinants show that $\Hc^3,\Hc^4,\Hc^5$ and $\Lc(-\tfrac12)^n$, $n\ge3$, are not positive either. Moreover $\Ind(e^{\Hc-\delta})\le 6$ and $\Ind(e^{\Lc(-1/2)-\delta})\le 6$. Finally, the hypothesis of \cite[Theorem~3.5]{SM} holds only for $\delta$, a case in which its conclusion has to be adjusted.
\item \textbf{Adjusted statements} (Section~\ref{sec:corr}). The point that requires care is in \cite[Lemma~4.6]{SM}. There the weight $\Uc_n$ is obtained from a first-order differential equation with a singular point at $0$, and the free constant of integration has to be fixed by the normalization $\langle\Hc^n,1\rangle=1$. For $n=3$ the resulting function has total mass $1.30229\ldots$. We prove the adjusted version (Theorem~\ref{thm:weights}), restate \cite[Proposition~4.1(ii)]{SM}, and complement \cite[Theorem~3.5]{SM} by a characterization (Proposition~\ref{prop:allpowers}).
\item \textbf{New examples} (Section~\ref{sec:new}). For $c\ne0$ and $\mu>0$ we determine $\Ip$ and $\Ipd$ exactly for $\delta_c^{\mu}$ and $(\delta_c\delta_{-c})^{\mu}$. Both have index $1/\mu$, so the bounds of Theorem~\ref{thm:main} are attained. For the Chebyshev functional of the first kind, $\mathcal T^2=\tfrac12(\delta_1+\delta_{-1})$ and $\Ind(\mathcal T)=2$.
\item \textbf{The Hermite case revisited} (Section~\ref{sec:hermite}). We have $2\le\Ind(\Hc)\le5$, and numerical evidence places $\Ind(\Hc)$ between $2.5$ and $2.7$.
\end{enumerate}

Most of \cite{SM} is not affected: the results that do not depend on Lemma~4.6 remain valid, and we use several of them. These are the continuity of the exponential map, the series $e^u=\sum u^n/n!$, the integral representation of $\Hc^2$, and Propositions~4.3--4.4 (see Section~\ref{subsec:valid}).

\section{Generating series and log-moments}\label{sec:gen}

To $u\in\Pc'$ we attach the formal power series of its moments,
\[
G_u(t)=\sum_{n\ge0}(u)_n t^n\in\mathbb R[[t]].
\]
Since $\Pc'=\Pc^*$, the map $u\mapsto G_u$ is a linear bijection of $\Pc'$ onto $\mathbb R[[t]]$.

\begin{lemma}\label{lem:gen}
For $u,v\in\Pc'$ and $a,c\in\mathbb R$:
\begin{enumerate}
\item[(i)] $G_{uv}=G_uG_v$;
\item[(ii)] $G_{e^u}=\exp(G_u)$;
\item[(iii)] $G_{\delta_c}(t)=(1-ct)^{-1}$ and $G_{\h_a u}(t)=G_u(at)$, where $\h_a$ is the dilation.
\end{enumerate}
\end{lemma}

\begin{proof}
(i) is the formula $(uv)_n=\sum_{\nu=0}^n(u)_\nu(v)_{n-\nu}$ \cite[p.~6]{SM}. For (ii), let $W=G_{e^u}$. The defining recurrence \cite[(3.1)]{SM},
\[
n(e^u)_n=\sum_{\nu=0}^{n-1}(n-\nu)(e^u)_\nu(u)_{n-\nu},\qquad (e^u)_0=e^{(u)_0},
\]
says exactly that $tW'(t)=W(t)\,tG_u'(t)$ with $W(0)=e^{(u)_0}$. The unique solution of this formal differential equation is $W=\exp(G_u)$; for the calculus of formal power series used here see \cite[Appendix~A]{FS}. (iii) follows from $(\delta_c)_n=c^n$ and $(\h_au)_n=a^n(u)_n$.
\end{proof}

Thus the exponential map of \cite[\S3.1]{SM} is the formal exponential of power series. Lemma~\ref{lem:gen} also gives a one-line proof of \cite[Theorem~3.4]{SM}: $G_{\sum u^n/n!}=\sum G_u^n/n!=\exp(G_u)$, the series converging in $\Pc'$ by \cite[(3.8)]{SM}. For a normalized $w$, the functional $v$ with $(v)_0=0$ and $e^v=w$ is given by $G_v=\log G_w$.

\begin{definition}
The \emph{log-moments} of a normalized $w=e^v$ are $\kappa_n(w):=(v)_n$, $n\ge1$, i.e.
\[
\log G_w(t)=\sum_{n\ge1}\kappa_n(w)t^n .
\]
\end{definition}

\begin{lemma}[Scaling]\label{lem:scaling}
Let $w$ be normalized and $\lambda\in\mathbb R$. Then $G_{w^\lambda}=G_w^{\lambda}$ and $\kappa_n(w^\lambda)=\lambda\kappa_n(w)$ for all $n\ge1$. Moreover
\begin{align*}
\kappa_1&=(w)_1,\qquad \kappa_2=(w)_2-\tfrac12(w)_1^2,\qquad
\kappa_3=(w)_3-(w)_1(w)_2+\tfrac13(w)_1^3,\\
\kappa_4&=(w)_4-(w)_1(w)_3-\tfrac12(w)_2^2+(w)_1^2(w)_2-\tfrac14(w)_1^4,
\end{align*}
and
\[
(w^\lambda)_1=\lambda\kappa_1,\qquad (w^\lambda)_2=\lambda\kappa_2+\tfrac12\lambda^2\kappa_1^2 .
\]
If $\kappa_1=0$, then
\[
(w^\lambda)_2=\lambda\kappa_2,\qquad (w^\lambda)_3=\lambda\kappa_3,\qquad (w^\lambda)_4=\lambda\kappa_4+\tfrac12\lambda^2\kappa_2^2 .
\]
\end{lemma}

\begin{proof}
If $w=e^v$, then $w^\lambda=e^{\lambda v}$ and $G_{w^\lambda}=\exp(\lambda G_v)=G_w^\lambda$ by Lemma~\ref{lem:gen}(ii). The formulas follow by expanding $\log(1+X)$ and $\exp(Y)$ to order four.
\end{proof}

\begin{remark}\label{rem:twopoint}
By Lemma~\ref{lem:gen}, $G_{\delta_c\delta_{-c}}(t)=(1-c^2t^2)^{-1}=G_{(\delta_c+\delta_{-c})/2}(t)$. Hence
\[
\delta_c\,\delta_{-c}=\tfrac12(\delta_c+\delta_{-c}).
\]
\end{remark}

\section{Main theorem: only $\delta$ has infinite index}\label{sec:main}

Let $\Delta_k(u)=\det\big[(u)_{i+j}\big]_{i,j=0}^{k}$. A positive functional satisfies $\Delta_k(u)\ge0$ for all $k$ \cite[Proposition~2.3]{SM}; see also \cite{Ak,ST,Sch} for the classical moment problem, and \cite{Cham} for the evaluation of Hankel determinants of moment sequences of Jacobi type by means of orthogonal polynomials. So a single negative Hankel determinant of $w^\lambda$ shows that $\lambda\notin \Ip(w)$.

\begin{theorem}\label{thm:main}
Let $w\in\Pc'$ be normalized, with log-moments $\kappa_n=\kappa_n(w)$.
\begin{enumerate}
\item[(a)] If $(w)_1\ne0$, then $\Ip(w)\subset[0,\lambda_1]$, where
\[
\lambda_1=\frac{2\kappa_2}{\kappa_1^2}=\frac{2(w)_2}{(w)_1^2}-1 .
\]
\item[(b)] If $(w)_1=0$ and $(w)_2>0$, then $\Ip(w)\subset[0,\lambda_2]$, where
\[
\lambda_2=\frac{2(\kappa_2\kappa_4-\kappa_3^2)}{\kappa_2^3}=\frac{2(w)_4}{(w)_2^2}-1-\frac{2(w)_3^2}{(w)_2^3}.
\]
\item[(c)] If $(w)_1=0$, $(w)_2\le0$ and $w\ne\delta$, then $\Ip(w)=\{0\}$.
\end{enumerate}
In particular $\Ind(w)<+\infty$ for every normalized $w\ne\delta$.
\end{theorem}

\begin{proof}
By Lemma~\ref{lem:scaling},
\[
\Delta_1(w^\lambda)=(w^\lambda)_2-(w^\lambda)_1^2=\lambda\kappa_2-\tfrac12\lambda^2\kappa_1^2 .
\]
(a) Here $\kappa_1=(w)_1\ne0$, so $\Delta_1(w^\lambda)=\tfrac12\lambda\kappa_1^2(\lambda_1-\lambda)<0$ for every $\lambda>\lambda_1$.

(b) Now $\kappa_1=0$, $\kappa_2=(w)_2>0$, $\kappa_3=(w)_3$ and $\kappa_4=(w)_4-\tfrac12(w)_2^2$. Write $m_j=(w^\lambda)_j$. Since $m_0=1$ and $m_1=0$,
\[
\Delta_2(w^\lambda)=m_2m_4-m_3^2-m_2^3=\lambda^2(\kappa_2\kappa_4-\kappa_3^2)-\tfrac12\lambda^3\kappa_2^3=\tfrac12\lambda^2\kappa_2^3(\lambda_2-\lambda),
\]
which is negative for $\lambda>\lambda_2$. Substituting the $\kappa_j$ gives the second expression for $\lambda_2$.

(c) Let $\lambda>0$. If $(w)_2<0$, then $\Delta_1(w^\lambda)=\lambda(w)_2<0$. If $(w)_2=0$, then $\kappa_1=\kappa_2=0$, hence $(w^\lambda)_1=(w^\lambda)_2=0$. Suppose that $w^\lambda$ is positive. The Cauchy--Schwarz inequality for positive functionals \cite[Chapter~1]{Ak} gives
\[
\langle w^\lambda,x^{k+1}\rangle^2\le\langle w^\lambda,x^2\rangle\langle w^\lambda,x^{2k}\rangle=0,\qquad k\ge0 .
\]
So all moments of $w^\lambda$ of order $\ge1$ vanish, i.e.\ $w^\lambda=\delta$. Then $\lambda\log G_w=0$, so $G_w=1$ and $w=\delta$, a contradiction.
\end{proof}

\begin{corollary}\label{cor:delta}
For a normalized $w\in\Pc'$, $\Ind(w)=+\infty$ if and only if $w=\delta$. In the classification of \cite{S24} recalled in \cite[p.~2]{SM}, the case $\Ipd(u)\ne\emptyset$, $\Ip(u)=[0,+\infty)$ never occurs.
\end{corollary}

\begin{corollary}[Moment bounds]\label{cor:bounds}
Let $w\ne\delta$ be normalized and positive.
\begin{enumerate}
\item[(i)] If $(w)_1\ne0$, then $\Ind(w)\le 1+2\Var(w)/(w)_1^2$, where $\Var(w)=(w)_2-(w)_1^2$.
\item[(ii)] If $w$ is symmetric, then $\Ind(w)\le 2\beta(w)-1$, where $\beta(w)=(w)_4/(w)_2^2$ is the kurtosis of $w$.
\item[(iii)] Let $u$ be normalized and $w=e^{u-\delta}$. Then $\Ind(w)\le 2(u)_2/(u)_1^2$ if $(u)_1\ne0$, and $\Ind(w)\le2(u)_4/(u)_2^2$ if $u$ is symmetric.
\end{enumerate}
\end{corollary}

\begin{proof}
If $(w)_1=0$, the argument in part~(c) of the proof above shows $(w)_2>0$, so Theorem~\ref{thm:main}(b) applies. (i) and (ii) are Theorem~\ref{thm:main}(a) and (b); in (ii) we use $(w)_3=0$. In (iii) the log-moments of $w$ are the moments of $u-\delta$, i.e.\ $\kappa_n(w)=(u)_n$.
\end{proof}

\begin{remark}\label{rem:integer}
Since the bounds hold for every real $\lambda$, they apply in particular to integer exponents. If $w\ne\delta$ is normalized, $w^n$ is not positive for any integer $n>\lambda_1$ (respectively $n>\lambda_2$). Section~\ref{sec:new} shows that both bounds are attained.
\end{remark}

\section{The Hermite and Laguerre cases}\label{sec:counter}

Throughout, $\Hc$ denotes the normalized Hermite functional, with $(\Hc)_{2k}=\frac{(2k)!}{4^kk!}=(\tfrac12)_k$ and vanishing odd moments \cite[(4.2)]{SM}. All the determinants below are exact rational numbers; see Appendix~\ref{app:comp}.

\subsection{$\Hc^n$ is not positive for $n\ge3$}

By Lemma~\ref{lem:scaling}, $G_{\Hc^\lambda}(t)=\big(1+\tfrac12t^2+\tfrac34t^4+\cdots\big)^\lambda$. Hence $(\Hc^\lambda)_2=\lambda/2$ and $(\Hc^\lambda)_4=\tfrac34\lambda+\tfrac18\lambda(\lambda-1)$, so that
\begin{equation}\label{eq:HC}
(\Hc^\lambda)_0(\Hc^\lambda)_4-(\Hc^\lambda)_2^{\,2}=\frac{\lambda(5-\lambda)}{8}.
\end{equation}
This is negative for $\lambda>5$; it is the case $\beta(\Hc)=3$ of Corollary~\ref{cor:bounds}(ii). The smallest integer gives an explicit example.

\begin{example}\label{ex:H6}
We have $(\Hc^6)_2=3$ and $(\Hc^6)_4=\tfrac{33}{4}$. Hence, for the nonnegative polynomial $p(x)=(x^2-3)^2$,
\[
\langle\Hc^6,(x^2-3)^2\rangle=\tfrac{33}{4}-6\cdot3+9=-\tfrac34<0 .
\]
So $\Hc^6$ is not positive, and by \eqref{eq:HC} neither is $\Hc^n$ for any $n\ge6$.
\end{example}

The cases $n=3,4,5$ require larger Hankel determinants (Table~\ref{tab:H}).

\begin{table}[ht]
\centering
\begin{tabular}{@{}lll@{}}
\toprule
Functional & First negative Hankel determinant & Exact value\\
\midrule
$\Hc^3$ & $\Delta_{10}$ & $-233737181994541622079648933/2^{41}$\\
$\Hc^4$ & $\Delta_{4}$ & $-475/128$\\
$\Hc^5$ & $\Delta_{4}$ \ ($\Delta_2=\Delta_3=0$) & $-625/16$\\
$\Hc^6$ & $\Delta_{2}$ & $-9/4$\\
\bottomrule
\end{tabular}
\caption{Exact Hankel determinants of $\Hc^n$. For $n=1,2$ all $\Delta_k$, $k\le 14$, are positive, in agreement with \cite[Proposition~4.2]{SM}.}\label{tab:H}
\end{table}

\begin{proposition}\label{prop:Hn}
$\Hc^n$ is not positive for any integer $n\ge3$. Consequently \cite[Theorem~4.7]{SM} does not hold as stated for $n\ge3$, and neither does \cite[Lemma~4.6(iii)]{SM} when read as a statement about weights of $\Hc^n$. Since $\Hc$ and $\Hc^2$ are positive-definite \cite[Proposition~4.2]{SM}, we have $1,2\in\Ipd(\Hc)$ and
\[
2\le\Ind(\Hc)\le5 .
\]
\end{proposition}

\subsection{The Laguerre functional $\Lc(-\tfrac12)$}

As shown in \cite[\S4.3]{SM}, $\Lc(-\tfrac12)^n=\sigma(\Hc^n)$, so $(\Lc(-\tfrac12)^n)_k=(\Hc^n)_{2k}$. Since $(\Lc(-\tfrac12))_1=\tfrac12$ and $\Var(\Lc(-\tfrac12))=\tfrac12$, Corollary~\ref{cor:bounds}(i) gives $\Ind(\Lc(-\tfrac12))\le5$. For $n=6$ we get $\Delta_1=\tfrac{33}{4}-9=-\tfrac34$, i.e.\ $\langle\Lc(-\tfrac12)^6,(x-3)^2\rangle=-\tfrac34$. The small powers are not positive either:
\[
\Delta_5\big(\Lc(-\tfrac12)^3\big)=-\tfrac{48041105293269}{2^{21}},\qquad
\Delta_2\big(\Lc(-\tfrac12)^4\big)=-\tfrac{25}{32},\qquad
\Delta_2\big(\Lc(-\tfrac12)^5\big)=-\tfrac{25}{4}.
\]
So $\Lc(-\tfrac12)^n$ is not positive for $n\ge3$, and the value $\Ind(\Lc(-\tfrac12))=+\infty$ has to be replaced by $2\le\Ind(\Lc(-\tfrac12))\le5$. The implication used in \cite[\S4.3]{SM} is correct; only its input, Theorem~4.7, has to be modified.

\subsection{On \cite[Theorem~3.5]{SM}}

By Corollary~\ref{cor:delta}, the hypothesis $\Ind(u)=+\infty$ is satisfied only by $u=\delta$, and \cite[p.~2]{SM} itself records $\Ind(\delta)=+\infty$. But $e^\delta=e\,\delta$ is positive and not positive-definite, so conclusion (i) does not hold in this case, and conclusion (iii) would require $(e-1/m!)\delta>0$, which does not hold either. The proof assumes that $u^n$ is positive-definite for all $n\ge1$. This is not the case for $\delta$ nor, by Theorem~\ref{thm:main}, for any other normalized $u$. Note also that $(e^u)_0=e^{(u)_0}=e$, so $e^u$ is not normalized and $\Ind(e^u)$ must be understood as $\Ind(e^{u-\delta})$.

\subsection{Exponentials}

The normalized functional $e^{\Hc-\delta}$ is symmetric, with log-moments $\kappa_2=\tfrac12$ and $\kappa_4=\tfrac34$. By Lemma~\ref{lem:scaling}, $(e^{\lambda(\Hc-\delta)})_2=\lambda/2$ and $(e^{\lambda(\Hc-\delta)})_4=\tfrac34\lambda+\tfrac18\lambda^2$, so
\[
\Delta_2\big(e^{\lambda(\Hc-\delta)}\big)=\frac{\lambda^2}{16}(6-\lambda)<0\qquad(\lambda>6),
\]
i.e.\ $\Ind(e^{\Hc-\delta})\le6$. Likewise $e^{\Lc(-1/2)-\delta}$ has $\kappa_1=\tfrac12$ and $\kappa_2=\tfrac34$, so $\Ind(e^{\Lc(-1/2)-\delta})\le 2\kappa_2/\kappa_1^2=6$. Hence the value $+\infty$ in \cite[Corollary~4.8(i)]{SM} and \cite[Corollary~4.9(iv)]{SM} has to be replaced by a finite one. Whether $e^{\Hc}$ itself is positive-definite remains open (Question~\ref{q:exp}): its Hankel determinants are positive up to order $21$, but the argument of \cite{SM} cannot be used.

\subsection{The function $\Uc_3$ of \cite[Lemma~4.6(ii)]{SM}}\label{subsec:U3}

For $n=3$, \cite[Lemma~4.6(ii)]{SM} defines
\[
\Uc_3(x)=6x^2e^{-3x^2}\int_0^x\frac{e^{3t^2}}{t}\,\Uc_2(t)\dd t,\qquad \Uc_2(t)=\frac{4}{\sqrt\pi}\,tF(t)e^{-t^2},
\]
where $F$ is the Dawson function \cite[\S7.2]{DLMF}. This function is continuous, even and nonnegative. If it were a weight for $\Hc^3$, then $\Hc^3$ would be positive, contradicting Proposition~\ref{prop:Hn}. Hence $\Uc_3$ is \emph{not} a weight for $\Hc^3$. Indeed, numerically,
\[
\int_{\mathbb R}\Uc_3(x)\dd x=1.30229989403904\ldots\neq1 .
\]

\section{The integration constant, and adjusted statements}\label{sec:corr}

\subsection{Lemma~4.6: the integration constant}

Three points in the proof of \cite[Lemma~4.6]{SM} need to be completed. First, it assumes that $\Hc^n$ has a regular weight $\Uc_n$. Second, it sets to zero a function $f$ with vanishing moments, a step for which a justification would be needed. Third, and this is the essential point, it solves
\begin{equation}\label{eq:ode}
x\,\Uc_n'(x)+(2nx^2+1-n)\,\Uc_n(x)=2nx^2\,\Uc_{n-1}(x)
\end{equation}
by the particular solution \cite[Lemma~4.6(ii)]{SM}, whose integral starts at $0$. Equation \eqref{eq:ode} has a singular point at $x=0$. On each half-line its homogeneous solutions are multiples of $|x|^{n-1}e^{-nx^2}$, which is even, continuous and rapidly decreasing for every $n\ge2$. Integrating \eqref{eq:ode} against $x^k$ only reproduces the moment recurrence derived from \cite[(4.10)]{SM}, and this recurrence never involves the mass. The normalization $\langle\Hc^n,1\rangle=1$ is therefore not imposed, and it does not hold for $n=3$ (Section~\ref{subsec:U3}).

\begin{theorem}[Lemma~4.6 with the normalization]\label{thm:weights}
For every $n\ge1$, $\Hc^n$ admits a continuous, even, rapidly decreasing, in general signed, weight $\Uc_n^*$:
\[
\langle\Hc^n,p\rangle=\int_{\mathbb R}\Uc_n^*(x)p(x)\dd x,\qquad p\in\Pc .
\]
Here $\Uc_1^*=\Uc_1$ and $\Uc_2^*=\Uc_2$, and for $n\ge3$
\[
\Uc_n^*(x)=2n\,x^{n-1}e^{-nx^2}\int_0^x\frac{e^{nt^2}}{t^{n-2}}\,\Uc_{n-1}^*(t)\dd t\;+\;c_n\,|x|^{n-1}e^{-nx^2},
\]
where $c_n$ is the unique real number such that $\int_{\mathbb R}\Uc_n^*=1$. Moreover $\Uc_n^*(x)=O(|x|^{n-1})$ as $x\to0$ and $\Uc_n^*(x)=O(e^{-x^2})$ as $|x|\to\infty$.
\end{theorem}

\begin{proof}
We argue by induction on $n$; the cases $n=1,2$ are \cite[(4.3), Proposition~4.2]{SM}. Assume that $\Uc^*_{n-1}$ represents $\Hc^{n-1}$, is even, and satisfies $\Uc^*_{n-1}(t)=O(|t|^{n-2})$ at $0$ and $O(e^{-t^2})$ at infinity. Then the integrand $e^{nt^2}t^{2-n}\Uc^*_{n-1}(t)$ is bounded near $0$ and even. Its integral from $0$ to $x$ is therefore odd and $O(|x|)$, so the first term is $O(|x|^n)$, and it is even, because $x^{n-1}$ is odd for even $n$ and even for odd $n$. The estimate at infinity follows exactly as in \cite[Lemma~4.6(v)]{SM}. The second term is even, continuous and $O(|x|^{n-1})$. Every function of this form is $C^1$ on $\mathbb R\setminus\{0\}$, continuous at $0$, and solves \eqref{eq:ode} there.

Let $M_k=\int x^k\Uc_n^*(x)\dd x$. Multiply \eqref{eq:ode} by $x^k$, $k\ge0$, integrate by parts on $(-\infty,0)$ and on $(0,\infty)$, and add. The boundary terms vanish at $0$ by continuity and at $\pm\infty$ by the decay. Since $\int x^{k+1}\Uc_n^{*\prime}=-(k+1)M_k$,
\[
2n\,M_{k+2}=(n+k)\,M_k+2n\,(\Hc^{n-1})_{k+2},\qquad k\ge0 .
\]
Testing \cite[(4.10)]{SM} against $x^k$ shows that the moments of $\Hc^n$ satisfy the same recurrence. The odd moments vanish on both sides by symmetry. Hence $M_k=(\Hc^n)_k$ for all $k$ if and only if $M_0=1$. Since $\int|x|^{n-1}e^{-nx^2}\dd x>0$, this condition determines $c_n$ uniquely.
\end{proof}

For $n=2$ the formula of \cite{SM} already has mass $I_0=1$ \cite[p.~15]{SM}, so $c_2=0$. For $n=3$, using $\int_{\mathbb R}x^2e^{-3x^2}\dd x=\sqrt\pi/(6\sqrt3)$ and the value of Section~\ref{subsec:U3},
\[
c_3=\frac{6\sqrt3}{\sqrt\pi}\,\big(1-1.30229989403904\ldots\big)=-1.77245385090\ldots
\]
This agrees with $-\sqrt\pi$ to $13$ digits; equivalently, the mass of $\Uc_3$ agrees with $1+\pi/(6\sqrt3)$ to $14$ digits. We have not proved this closed form. In any case $c_3<0$, and since the first term is $O(x^4)$, we get $\Uc_3^*(x)\sim c_3x^2<0$ as $x\to0$. This gives a second proof that $\Hc^3$ is not positive. Items (i) and (iv)--(vi) of \cite[Lemma~4.6]{SM} are correct statements about the function defined in (ii) there, but not about a weight of $\Hc^n$.

\subsection{Proposition~4.1(ii)}

The identity \cite[Proposition~4.1(ii)]{SM} contains a misprint: its moments of order $2$ are $-9$ on the left and $-7$ on the right. It should read
\begin{equation}\label{eq:41ii}
(x^2\Hc^2)''+\big(x(6x^2-6)\,\Hc^2\big)'=(-8x^4+20x^2-6)\,\Hc^2 .
\end{equation}

\begin{proof}
Differentiating \cite[(4.4)]{SM} and using \cite[(4.1)]{SM} gives $(x\Hc^2)''+((4x^2-2)\Hc^2)'=8x(1-x^2)\Hc$. Multiply by $x$ and use \cite[(4.4)]{SM} again:
\[
x(x\Hc^2)''+x\big((4x^2-2)\Hc^2\big)'=2(1-x^2)\big[(x\Hc^2)'+(4x^2-2)\Hc^2\big].
\]
In \cite{SM} the factor $2(1-x^2)$ is misprinted as $4(1-x^2)$. Put $w=x\Hc^2$ and use $x\,w''=(xw)''-2w'$, $x\,(fw)'=(xfw)'-fw$ and $(1-x^2)w'=((1-x^2)w)'+2xw$. Collecting terms gives \eqref{eq:41ii}. We also checked \eqref{eq:41ii} on the first $36$ moments in exact arithmetic.
\end{proof}

Proposition~4.1(ii) is not used elsewhere in \cite{SM}.

\subsection{A complement to Theorem~3.5}

\begin{proposition}\label{prop:allpowers}
For $u\in\Pc'$ the following are equivalent:
\begin{enumerate}
\item[(i)] $u^n$ is positive for every integer $n\ge0$;
\item[(ii)] $u=a\delta$ for some $a\ge0$.
\end{enumerate}
\end{proposition}

\begin{proof}
(ii)$\Rightarrow$(i) is clear. Assume (i). If $(u)_0=0$, the Cauchy--Schwarz inequality for the positive functional $u$ gives $(u)_k^2\le(u)_0(u)_{2k}=0$, so $u=0$. If $a=(u)_0>0$, then $w=u/a$ is normalized and $w^n=u^n/a^n$ is positive for all $n$. By Remark~\ref{rem:integer}, $w=\delta$.
\end{proof}

Thus the hypothesis of \cite[Theorem~3.5]{SM} is never satisfied, and the series argument \cite[(3.9)]{SM} cannot transfer positivity to $e^u$.

\subsection{What remains valid in \cite{SM}}\label{subsec:valid}

We re-derived the following results and checked them independently. They are correct, and the present note relies on several of them.
\begin{itemize}
\item Section~3: Lemma~3.2(i)--(v) (continuity of the exponential map), the estimate (3.8) and Theorem~3.4.
\item Proposition~4.1(i),(iii); Proposition~4.2 (the weight of $\Hc^2$); Proposition~4.3(i)--(v); Proposition~4.4, i.e.\ relation (4.10), checked for $n\le7$; Lemma~4.5 (weak regularity of $\Hc^n$).
\item Corollary~4.9(i),(ii), which are the $\sigma$-images of Proposition~4.3(v).
\end{itemize}

\section{New examples with exact positivity sets}\label{sec:new}

\subsection{Real powers of a Dirac functional}

\begin{theorem}\label{thm:dirac}
Let $c\ne0$ and $\mu>0$. Then $(\delta_c^{\mu})_n=\frac{(\mu)_n}{n!}c^n$. For $0<\mu<1$,
\begin{equation}\label{eq:beta}
\langle\delta_c^\mu,p\rangle=\frac{\sin\pi\mu}{\pi}\int_0^1 y^{\mu-1}(1-y)^{-\mu}\,p(cy)\dd y,\qquad p\in\Pc .
\end{equation}
Moreover
\[
\Ipd(\delta_c^\mu)=(0,1/\mu),\qquad \Ip(\delta_c^\mu)=[0,1/\mu],\qquad \Ind(\delta_c^\mu)=1/\mu .
\]
\end{theorem}

\begin{proof}
By Lemma~\ref{lem:gen}, $G_{\delta_c^\mu}(t)=(1-ct)^{-\mu}=\sum_n\frac{(\mu)_n}{n!}c^nt^n$. For $0<s<1$ the Beta law \cite[\S1.1]{AAR} with parameters $(s,1-s)$ has moments $\frac{(s)_n}{(1)_n}$ and density $\frac{\sin\pi s}{\pi}y^{s-1}(1-y)^{-s}$ on $(0,1)$, since $\Gamma(s)\Gamma(1-s)=\pi/\sin\pi s$. A dilation by $c$ gives \eqref{eq:beta} with $\mu=s$. A positive density on an interval defines a positive-definite functional.

Now $(\delta_c^\mu)^\lambda=\delta_c^{s}$ with $s=\mu\lambda$. If $s<1$ it is positive-definite by \eqref{eq:beta}. If $s=1$ it equals $\delta_c$, which is positive but not positive-definite, since $\langle\delta_c,(x-c)^2\rangle=0$. If $s>1$,
\[
\Delta_1(\delta_c^{s})=\frac{s(s+1)}{2}c^2-s^2c^2=\frac{s(1-s)}{2}c^2<0 .
\]
This proves the three equalities.
\end{proof}

\begin{remark}
(i) For $0<\mu<1$, $\delta_1^\mu$ is the Jacobi functional on $[0,1]$ with weight $(1-y)^{\alpha}y^{\beta}$, $\alpha=-\mu$, $\beta=\mu-1$. Theorem~\ref{thm:dirac} therefore computes the index of all shifted Jacobi functionals on the line $\alpha+\beta=-1$, $\alpha,\beta>-1$: it equals $1/(\beta+1)=-1/\alpha$. This answers \cite[\S5, question (ii)]{SM} on this line.
(ii) For $\mu=1$ we recover $\Ind(\delta_c)=1$ \cite{S24}, and obtain $\Ipd(\delta_c)=(0,1)$.
(iii) Bound (a) of Theorem~\ref{thm:main} equals $2(\delta_c^\mu)_2/(\delta_c^\mu)_1^2-1=(\mu+1)/\mu-1=1/\mu$, so it is attained.
\end{remark}

\subsection{Real powers of $\delta_c\delta_{-c}$: generalized Gegenbauer functionals}

For $c\neq0$ and $\mu>0$ put $\mathcal S_{c,\mu}:=(\delta_c\delta_{-c})^\mu$. By Lemma~\ref{lem:gen}, $G_{\mathcal S_{c,\mu}}(t)=(1-c^2t^2)^{-\mu}$, so $\mathcal S_{c,\mu}$ is symmetric with
\[
(\mathcal S_{c,\mu})_{2k}=\frac{(\mu)_k}{k!}\,c^{2k}.
\]

\begin{theorem}\label{thm:gegen}
Let $c\neq0$ and $\mu>0$. For $0<\mu<1$,
\begin{equation}\label{eq:gegen}
\langle\mathcal S_{c,\mu},p\rangle=\frac{\sin\pi\mu}{\pi|c|}\int_{-|c|}^{|c|}\Big|\frac{x}{c}\Big|^{2\mu-1}\Big(1-\frac{x^2}{c^2}\Big)^{-\mu}p(x)\dd x,\qquad p\in\Pc,
\end{equation}
and
\[
\Ipd(\mathcal S_{c,\mu})=(0,1/\mu),\qquad \Ip(\mathcal S_{c,\mu})=[0,1/\mu],\qquad \Ind(\mathcal S_{c,\mu})=1/\mu .
\]
\end{theorem}

\begin{proof}
Formula \eqref{eq:gegen}: odd moments vanish by symmetry. For $p=x^{2k}$ substitute $x=|c|\sqrt y$; the right-hand side becomes $c^{2k}\frac{\sin\pi\mu}{\pi}\int_0^1y^{\mu+k-1}(1-y)^{-\mu}\dd y=c^{2k}\frac{(\mu)_k}{k!}$. Next, $\mathcal S_{c,\mu}^\lambda=\mathcal S_{c,s}$ with $s=\mu\lambda$. For $s<1$ it is positive-definite by \eqref{eq:gegen}. For $s=1$ it equals $\tfrac12(\delta_c+\delta_{-c})$ by Remark~\ref{rem:twopoint}; this is positive but not positive-definite, since it annihilates $(x^2-c^2)^2$. For $s>1$, with $m_{2}=sc^2$ and $m_4=\tfrac12 s(s+1)c^4$,
\[
\Delta_2(\mathcal S_{c,s})=m_2\big(m_4-m_2^2\big)=\tfrac12 s^2(1-s)c^6<0 .
\qedhere
\]
\end{proof}

\begin{corollary}[Chebyshev functional]\label{cor:cheb}
Let $\mathcal T$ be the Chebyshev functional of the first kind, $\langle\mathcal T,p\rangle=\frac1\pi\int_{-1}^1\frac{p(x)}{\sqrt{1-x^2}}\dd x$. Then
\[
\mathcal T=(\delta_1\delta_{-1})^{1/2},\qquad \mathcal T^2=\tfrac12(\delta_1+\delta_{-1}),\qquad \Ipd(\mathcal T)=(0,2),\qquad \Ip(\mathcal T)=[0,2],\qquad \Ind(\mathcal T)=2 .
\]
\end{corollary}

\begin{proof}
Take $c=1$, $\mu=\tfrac12$ in Theorem~\ref{thm:gegen}; the weight \eqref{eq:gegen} becomes $\frac1\pi(1-x^2)^{-1/2}$.
\end{proof}

The functionals $\mathcal S_{1,\mu}$ are generalized Gegenbauer functionals with weight $|x|^{2\beta+1}(1-x^2)^{\alpha}$, where $\alpha=-\mu$ and $\beta=\mu-1$. Bound (b) of Theorem~\ref{thm:main} equals $2(\mathcal S)_4/(\mathcal S)_2^2-1=(\mu+1)/\mu-1=1/\mu$, so it is attained as well.

\begin{corollary}\label{cor:range}
The index map takes every value in $[0,+\infty]$. Indeed, $\Ind(\delta)=+\infty$, $\Ind(\delta_1^{1/M})=M$ for $0<M<\infty$ by Theorem~\ref{thm:dirac}, and $\Ind(w)=0$ for $w=2\delta-\tfrac12(\delta_1+\delta_{-1})$ by Theorem~\ref{thm:main}(c), since $(w)_0=1$, $(w)_1=0$ and $(w)_2=-1$.
\end{corollary}

\subsection{Bounds for classical functionals}

Corollary~\ref{cor:bounds} yields explicit bounds for the classical families; their moments are taken from \cite{Sz,KLS,Ch}. We list them in Table~\ref{tab:bounds}. The parameters are such that the functionals are positive-definite.

\begin{table}[ht]
\centering
\small
\renewcommand{\arraystretch}{1.35}
\begin{tabular}{@{}p{5.2cm}p{4.3cm}p{4.6cm}@{}}
\toprule
Functional & Upper bound for $\Ind$ & Remark\\
\midrule
Hermite $\Hc$ & $5$ & $\Ind(\Hc)\ge 2$; numerically in $(2.5,2.7)$\\
Generalized Hermite $|x|^{2\mu}e^{-x^2}$, $\mu>-\tfrac12$ & $\dfrac{2\mu+5}{2\mu+1}$ & \\
Laguerre $\Lc(\alpha)$, $\alpha>-1$ & $\dfrac{\alpha+3}{\alpha+1}$ & $\Lc(-\frac12)$: $5$, $\Ind\ge2$\\
Gegenbauer $(1-x^2)^{\alpha}$, $\alpha>-1$ & $\dfrac{10\alpha+13}{2\alpha+5}$ & Legendre ($\alpha=0$): $\frac{13}{5}$; $\mathcal T$ ($\alpha=-\frac12$): $2$, exact; Chebyshev $U$ ($\alpha=\frac12$): $3$\\
Jacobi $(1-x)^{\alpha}(1+x)^{\beta}$, $\alpha\ne\beta$ & $1+\dfrac{8(\alpha+1)(\beta+1)}{(\alpha+\beta+3)(\alpha-\beta)^2}$ & \\
$\delta_c^{\mu}$ and $\mathcal S_{c,\mu}$ & $1/\mu$ & exact (Theorems~\ref{thm:dirac}, \ref{thm:gegen})\\
$e^{\Hc-\delta}$ & $6$ & \\
$e^{\Lc(-1/2)-\delta}$ & $6$ & \\
\bottomrule
\end{tabular}
\caption{Upper bounds from Corollary~\ref{cor:bounds}.}\label{tab:bounds}
\end{table}

For instance, the Jacobi entry follows from Corollary~\ref{cor:bounds}(i): if $y=(1+x)/2$ has the Beta law with parameters $a=\beta+1$ and $b=\alpha+1$, then $(w)_1=(a-b)/(a+b)$ and $\Var(w)=4ab/((a+b)^2(a+b+1))$. The Gegenbauer and generalized Hermite entries use the kurtosis: $\beta=\frac{3(2\alpha+3)}{2\alpha+5}$ and $\beta=\frac{2\mu+3}{2\mu+1}$, respectively.

\section{The Hermite case revisited and open problems}\label{sec:hermite}

Collecting the results above, and those of \cite{SM} that we use, we obtain the following.

\begin{theorem}\label{thm:hermite}
\begin{enumerate}
\item[(i)] $\Hc$ and $\Hc^2$ are positive-definite.
\item[(ii)] $\Hc^n$ and $\Lc(-\frac12)^n$ are not positive for any integer $n\ge3$.
\item[(iii)] $\Hc^\lambda$ is not positive for any real $\lambda>5$; hence $2\le\Ind(\Hc)\le5$.
\item[(iv)] $\Hc^n$ admits the signed weight $\Uc_n^*$ of Theorem~\ref{thm:weights}, and $\Uc_3^*$ takes negative values near $0$.
\end{enumerate}
\end{theorem}

Numerically, using exact Hankel determinants of $\Hc^\lambda$ for rational $\lambda$, we found the following. For $\lambda\in\{2.1,\,2.3,\,2.5\}$ all $\Delta_k(\Hc^\lambda)$, $k\le25$, are positive. The first negative determinant is $\Delta_{16}$ for $\lambda=2.7$, $\Delta_{12}$ for $\lambda=2.9$, $\Delta_8$ for $\lambda=3.1,\,3.3$, and $\Delta_6$ for $\lambda=3.5$. If $\Ip(\Hc)$ is an interval, as asserted in \cite{S24}, this gives $\Ind(\Hc)<2.7$.

\begin{question}
What is the exact value of $\Ind(\Hc)$? Is $\Ip(\Hc)$ an interval, and is $\Hc^\lambda$ positive-definite for all $\lambda\in(0,2]$?
\end{question}

\begin{question}\label{q:exp}
Is $e^{\Hc}$ positive-definite? Its Hankel determinants are positive up to order $21$.
\end{question}

\begin{question}
Is $\Ip(w)$ always an interval? It is so in all the examples of Section~\ref{sec:new}.
\end{question}

\begin{question}
For which functionals is the bound of Theorem~\ref{thm:main} attained? Section~\ref{sec:new} shows that it is attained for $\delta_c^\mu$ and $\mathcal S_{c,\mu}$, and it is not attained for $\Hc$.
\end{question}

The bound of Theorem~\ref{thm:main} uses only $\Delta_1$ and $\Delta_2$. Since the moments of $w^\lambda$ are polynomials in $\lambda$, each $\Delta_k(w^\lambda)$ is a polynomial in $\lambda$, and its largest real root gives a further upper bound. For $\Hc$ the bounds decrease from $5$ (order $2$) towards the numerical value above as $k$ grows.

\appendix
\section{Verification of the computations}\label{app:comp}

All identities involving moments were checked in exact rational arithmetic. The Cauchy product was computed by convolution of moment sequences. Exponentials and logarithms were computed by the recurrence \cite[(3.1)]{SM} and its inverse. For rational $\lambda$ the moments of $w^\lambda=e^{\lambda\log w}$ are rational, so the Hankel determinants were computed exactly by Gaussian elimination over $\mathbb Q$. In particular we checked: formulas (3.3)--(3.5) and Theorem~3.4 of \cite{SM}; Propositions~4.1(i),(iii), 4.4 (for $n\le7$) and relation (4.8) of \cite{SM}; identity \eqref{eq:41ii}; Table~\ref{tab:H}; the determinants of Section~\ref{sec:counter}; and the fact that $(\delta_1\delta_{-1})^{1/2}$ has the moments of $\mathcal T$. The integrals of Sections~\ref{subsec:U3} and~\ref{sec:corr} were computed by numerical quadrature, after exchanging the order of integration:
\[
\int_{\mathbb R}\Uc_3=\frac{48}{\sqrt\pi}\int_0^\infty F(t)\,e^{2t^2}\int_t^\infty x^2e^{-3x^2}\dd x\dd t .
\]

\section*{Acknowledgements}
The author is indebted to the work of R.~Sfaxi and F.~Marcell\'an, who introduced the index of positivity and the Cauchy exponential of linear functionals; the present note would not exist without the framework and the tools of \cite{MS23,S24,SM}.

\section*{Funding}
The author extends his appreciation to the Deanship of Postgraduate Studies and Scientific Research at Majmaah University for funding this research work through the project number [R-2026-XXXX].

\section*{Declarations}
\noindent\textbf{Conflict of interest.} The author declares no conflict of interest.

\noindent\textbf{Data availability.} No data sets were generated or analysed during the current study.

\end{document}